\documentstyle{article}
\begin{document}
\setlength{\baselineskip}{15pt}

\mbox{}

\begin{center}
{\large {\bf Analytical quantile solution for the 
S-distribution, random number generation and statistical data modeling}} \\
\mbox{}\\
\mbox{}\\
\noindent {\bf Benito Hern\'{a}ndez--Bermejo$^{(1)}$
\hspace{1.5cm} Albert Sorribas$^{(2,*)}$}
\end{center}
{\em 
\begin{description}

\item \mbox{}

\item{$^{(1)}$} Departamento de F\'{\i}sica Matem\'{a}tica y Fluidos,
Universidad Nacional de Educaci\'{o}n a Distancia. Senda del Rey S/N, 28040 
Madrid, Spain. \\
E-mail: bhernand@apphys.uned.es

\item \mbox{}

\item{$^{(2)}$} Departament de Ci\`{e}ncies M\`{e}diques B\`{a}siques, 
Universitat de Lleida, Av. Rovira Roure 44, 25198-Lleida, Spain. \\
E-mail: Albert.Sorrribas@cmb.udl.es
\end{description} }

\mbox{}

\mbox{}

\mbox{}

\noindent {\bf Running Title: } Analytical solution of the S-distribution.

\mbox{}

\mbox{}

\noindent $^*$ Corresponding author.

\mbox{}

\mbox{}

\noindent Accepted for publication in {\it Biometrical Journal}
\newtheorem{fig}{Figure}

\pagebreak

\begin{center}
{\bf Summary}
\end{center}

The selection of a specific statistical distribution as a model for
describing the population behavior of a given variable is seldom a simple
problem. One strategy consists in testing different distributions 
(normal, lognormal, Weibull, etc.), and selecting the one providing 
the best fit to the observed data and being the most parsimonious. 
Alternatively, one can make a choice based
on theoretical arguments and simply fit the corresponding parameters 
to the observed data. In either case, different distributions can give similar 
results and provide almost equivalent models for a given data set. 
Model selection can be more complicated when the goal is to describe a trend 
in the distribution of a given variable. In those cases, changes in shape and 
skewness are difficult to represent by a single distributional form. 
As an alternative to the use of complicated families of distributions as
models for data, the S-distribution [{\sc Voit, E.O. }(1992) Biom.J. 7:855-878]
provides a highly flexible mathematical form in which the density is defined
as a function of the cumulative. S-distributions can accurately approximate many known 
continuous and unimodal distributions, preserving the well known limit relationships between 
them. Besides representing well-known distributions, S-distributions provide an infinity of 
new possibilities that do not correspond with known classical distributions. Although 
the utility and performance of this general form has been clearly proved 
in different applications, its definition as a differential equation is a potential 
drawback for some problems. In this paper we obtain an analytical solution for the quantile 
equation that highly simplifies the use of S-distributions. We show the utility of 
this solution in different applications. After classifying the different 
qualitative behaviors of the S-distribution in parameter space, 
we show how to obtain different S-distributions that
accomplish specific constraints. One of the most interesting cases is the
possibility of obtaining distributions that acomplish $P(X \leq X_c)=0$. Then,
we demonstrate that the quantile solution facilitates the use of S-distributions
in Monte-Carlo experiments through the generation of random samples. Finally, we
show how to fit an S-distribution to actual data, so that the resulting 
distribution can be used as a statistical model for them.

\mbox{}
\mbox{}


\pagebreak
\section{Introduction}

Statistical data modeling requires selecting a distribution among all the different theoretical 
models. While in some cases there is a sound justification for selecting a specific model, in many 
others no clear reason exists for selecting a particular one. In practice, different distributions
can lead to similar fits for a given data set and the selection of a specific distribution
is a matter of taste ({\sc Sorribas} et al., 2000). Data analysis in this context is a sort of 
trial and error problem until a compromise is reached. However, minute changes in the data sets may change the conclusions
by forcing the selection of a completely different distribution as a basic model. This has been observed
for instance in studies on contamination concentration in fisheries when the concentration is analyzed as
a function of fish lenght ({\sc Balthis} et al., 1996; {\sc Voit} et al., 1995). The same problem has been 
observed in modeling the change in weight distribution associated with age in which 
no single distribution model could account for all kind of shapes and skeweness that can be 
observed in the collected data ({\sc Sorribas} et al., 2000). In these cases, an added difficulty is that
the distribution must be truncated at zero since only positive values can be observed. In this context, it is
well known that the use of distributions such as the normal is not justified and that 
other alternatives must be explored. This problem is specially significant in the case 
of studying clinical parameters and trying to derive the corresponding normality intervals. 

Data transformation is one solution for these problems. The well-known Box-Cox transformation of the random variable is a convenient one. Other possibilities are the use of kernel estimation of the density ({\sc Wand and Jones}, 1995) or methods for generating distributions such as the {\it maximum entropy method} ({\sc Ryu}, 1993; {\sc Wagner and Geyer}, 1995). However, it would interesting to explore if some parametric alternative, in the sense of a family of distributions, could be used. At this point, it may be useful to find a single distributional form covering the entire spectrum of observed distributional shapes. This may be achieved by a superfamily of distributions that would contain all the required functions as particular cases. Although some of these families exist (see for instance {\sc Johnson} and {\sc Kotz}, 1970) in general they are so unwieldy that they are of very limited use in practical data analysis. 

As an alternative, the S-distribution family provides a highly flexible model that greatly facilitates 
data representation in these situations and provides a general framework for accurately approximating 
most unimodal distributions ({\sc Voit}, 1992; {\sc Voit} and {\sc Yu}, 1994). Furthermore, parameter estimation 
and computation of distribution properties is straighforward ({\sc Voit}, 1992; {\sc Voit}, 2000). 

The S-distribution was defined by Voit (see {\sc Voit}, 1992, for a detailed discussion on the rationale
for this definition) in terms of a differential equation, in which the cumulative, $F$, is the dependent variable:

\begin{equation}
\label{sddef}
	 \frac{\mbox{d} F}{\mbox{d} X}  = 
	 \alpha ( F^g - F^h ) \; ,   \;\:\;\: F(X_0) = F_0 \; , \;\:\;\: 
	 \alpha > 0 \; , \;\: h > g
\end{equation}

The density, {\it pdf}: $\mbox{d} F / \mbox{d} X $, is thus defined as a function of the cumulative,
$cdf$: $F$, and the random variable does not appear explicitly. Computation of {\it cdf} is made by 
numerically integrating (\ref{sddef}) and quantiles can be computed by integrating the 
inverse S-distribution defined as ({\sc Voit}, 1992):

\begin{equation}
\label{inverseSD}
	 \frac{\mbox{d} X}{\mbox{d} F}  = 
	 \frac{1}{\alpha ( F^g - F^h )}	  \; ,   \;\:\;\: X(F_0) = X_0 \;  
\end{equation}

For simplicity, an S-distribution with real parameters $\alpha$, $g$ and $h$, 
and initial condition $F(X_0)=F_0$, will be indicated as: $S[F_0,X_0,\alpha,g,h]$. According 
to its definition, it is worth recalling that $X_0$ is a parameter of location, provided a
value of $F_0$, and that $g$ and $h$ are parameters that define the shape of the distribution.
The parameter $\alpha$ is inversely related with the spread of the distribution ({\sc Voit},1992).
In many applications, it is useful to choose $F_0=0.5$, so that $X_0$ is the median of the
distribution. We shall use the median as a value for $X_0$ in most examples. 

The utility of the S-distribution has been demonstrated in different ways. First, it was shown that
many known unimodal statistical distributions can be accurately represented as an S-distribution and that the 
relationships between them are preserved when they are approximated by 
S-distributions ({\sc Voit}, 1992; {\sc Voit} and {\sc Yu}, 1994). Second, the usefulness of 
S-distributions as models for actual data has been explored in the case of risk assessment 
and standard curves, showing that S-distributions provide an accurate model for the observed data
({\sc Balthis} et al., 1996; {\sc Sorribas} et al., 2000; {\sc Voit} et al., 1995). Furtheremore, 
it was shown that S-distributions can be used to model the dynamical change in the statistical distribution 
of a given quantity as a result of growth ({\sc Voit} and {\sc Sorribas}, 2000). Finally, 
a maximum likelihood estimator has been recently obtained for the S-distribution shape parameters
({\sc Voit}, 2000).

In all cases, most of the computations require numerical integration of (\ref{sddef}), since no closed, 
usable analytical solution had been obtained for this equation (see however {\sc Voit} and {\sc Savageau}, 1984).
In this paper we develop a complete and useful analytical solution for the S-distribution
equation. Using this solution we discuss the existence of some special parameter values
that lead to particular cases and demonstrate some important features of the resulting distributions. 
The analytical solution obtained is of the form of a quantile equation $X \equiv X(F)$. Far from 
being an inconvenient, this solution greatly facilitates the use of S-distributions in
Monte-Carlo simulations and data modeling as we shall see in the final part of the article.

\section{Analytical solution of the S-distribution}

The problem we aim to solve is:
\begin{equation}
\label{prob2}
	 \frac{\mbox{d} F}{\mbox{d} X} \equiv \dot{F} = 
	 \alpha ( F^g - F^h ) \; , \;\:\;\: 
	 0 < F < 1  \; , \;\:\;\: F(X_0) = F_0 \; , \;\:\;\: 
	 \alpha > 0 \; , \;\: h > g
\end{equation}
Recall that in (\ref{prob2}) $\alpha$, $g$ and $h$ are real constants. 
With the definition of the parameters as given in (\ref{prob2}) we obtain 
$\dot{F} >0$ and therefore $F(X)$ is a strictly monotonous 
function. For convenience, we define $\gamma = h-g >0$. Then we can write 
equation (\ref{prob2}) as:
\begin{equation}
\label{prob3}
	 \dot{F} = \alpha F^g (1 - F^{\gamma} ) 
\end{equation}

This equation is separable:
\begin{equation}
\label{sol1}
	 X = X_0 + \frac{1}{\alpha}
	 \int_{F_0}^{F} \frac{\mbox{d}\xi}{\xi^g (1 - \xi^{\gamma} )}  
\end{equation}
We thus obtain $X$ as a function of $F$. Now notice that the integral can 
be exactly solved in the domain of interest, since:
\begin{equation}    
	 \frac{1}{ 1 - \xi^{\gamma} } = \sum_{k=0}^{\infty} \xi^{k \gamma}
\end{equation}
This series is convergent because $0<\xi^{\gamma}<1$ in the integration 
domain. Therefore it is immediate to show ({\sc Churchill} and {\sc Brown}, 1984) that term by term 
integration is valid here:
\begin{equation} 
\label{intser}
	 \int_{F_0}^{F} \frac{\mbox{d}\xi}{\xi^g (1 - \xi^{\gamma} )} = 
	 \sum_{k=0}^{\infty} \int_{F_0}^{F}  \xi^{k \gamma - g} \mbox{d} \xi
\end{equation}
Recall that this equality holds for $0<F<1$. The problem is thus exactly 
solvable. Two different cases must be distinguished now in parameter space.

\subsection{The generic case}

Given $g$ and $h$ in (\ref{prob2}), we say that we are in the {\em generic\/} 
or {\em nondegenerate\/} case if there is no natural $k$ ($k=0,1,2 \ldots$) 
such that $k \gamma - g = -1$. This includes all S-distributions with $g<1$ and 
most with $g>1$, while the case $g=1$ is nongeneric, as it shall be explained 
in what follows. It is important to announce that, as we shall see in detail at 
the end of this subsection and in the next one, all S-distributions of practical 
interest belong to the generic case.

Thus, in the generic case no integral in the infinite series of integrals in 
(\ref{intser}) produces a logarithmic contribution. Term by term integration 
is then straighforward. Taking (\ref{sol1}) into account and defining 
$\lambda =1-g$ we find:
\begin{equation}
\label{gen1}
	 X = X_0 + \frac{\xi^{\lambda}}{\alpha \lambda} \left. \left( 1 + 
	 \sum_{k=1}^{\infty} \frac{\xi^{k \gamma}}{1+k \gamma / \lambda} \right)
	 \right| _{F_0}^{F}
\end{equation}
We now make use of Lerch's transcendent ({\sc Erd\'{e}lyi} et al., 1953), which is defined as:
\begin{equation}
\label{lerch}
	\Phi [z,s,v] = \sum _{n=0}^{\infty} \frac{z^n}{(v+n)^s} \; , \;\:\;\: 
	\mid z \mid < 1  \; , \;\:\;\: v \neq 0, -1, -2, \ldots
\end{equation}
It is then an immediate task to prove that (\ref{gen1}) becomes:
\begin{equation}
\label{gen2}
	 X = X_0 +
	  \frac{F^{\lambda}}{\alpha \lambda} \left( 1 + 
	  \frac{\lambda F^{\gamma}}{\gamma} \Phi [F^{\gamma},1,1 + 
	  \lambda / \gamma] \right) -
	  \frac{F_0^{\lambda}}{\alpha \lambda} \left( 1 + 
	  \frac{\lambda F_0^{\gamma}}{\gamma} \Phi [F_0^{\gamma},1,1 + 
	  \lambda / \gamma] \right) 
\end{equation}
As a technical but important detail of this derivation, note that the 
identification
\begin{equation}
\label{iden}
	  \sum_{k=1}^{\infty} \frac{\xi^{k \gamma}}{1+k \gamma / \lambda} = 
	  \frac{\lambda \xi^{\gamma}}{\gamma} \Phi [\xi^{\gamma},1,1 + 
	  \lambda / \gamma] 
\end{equation}
will be valid only, according to (\ref{lerch}), when:

\mbox{}

\noindent {\bf i)} $\mid z \mid < 1$, which is trivially true: now we have 
$z=\xi^{\gamma}$, where $0 < \mbox{min} \{ F_0 , F \} \leq \xi \leq 
\mbox{max} \{ F_0, F \} <1$ and $\gamma >0$.

\mbox{}

\noindent {\bf ii)} $v \neq 0,-1,-2, \ldots$. What we have is $v = 1+
\lambda / \gamma$. But the condition $v = 0,-1,-2, \ldots$ implies trivially
that there exists a $k \in \{ 1,2,3 \ldots \}$ such that $kh-(k+1)g=-1$. This  
is a degeneracy condition for every $k$, and therefore now is excluded by 
hypothesis.

\mbox{}

Consequently, (\ref{iden}) is always valid in this context and (\ref{gen2}) 
is the general solution of the generic case. Notice that we have found the 
solution in the inverse form $X \equiv X(F)$, instead of $F \equiv F(X)$. 
However, both are completely equivalent in this case, because we are dealing 
with a monotonic, one-to-one solution. 

At this point, it is interesting to elaborate on the asymptotic behaviour of 
the generic solution (\ref{gen1}) or (\ref{gen2}). Such behaviour is not 
always the same and depends on the parameter values. Of course, this analysis 
is possible now since we have the explicit solution of the problem. Four 
different cases are to be considered:

\subsubsection{Case I: $\;\: g>1 $}

Let us look at the left side of the distribution (i.e., the limit 
$X \rightarrow - \infty$ or, equivalently, $F \rightarrow 0^+$). From 
(\ref{gen1}) we have that the dominant behaviour is:
\begin{equation}
	 X (F \rightarrow 0^+) \; \sim \; \frac{F^{\lambda}}{\alpha \lambda} 
	 + \mbox{constant}
\end{equation}
Taking into account that now $\lambda <0$ we find:
\begin{equation}
	 \lim _{F \rightarrow 0^+} X(F) \; = \; - \infty
\end{equation}
Then the distribution has an infinite left tail for this parameter range. As 
we shall see in what is to follow, this property is exceptional in the 
generic case.

\subsubsection{Case II: $\;\: 0<g<1$}

Mathematically, this is a very interesting case. Since $0< \lambda <1$, we 
have:
\begin{equation}
\label{noleft1}
	 \lim _{F \rightarrow 0^+} X(F) \; = \; X_0 -
	  \frac{F_0^{\lambda}}{\alpha \lambda} \left( 1 + 
	  \frac{\lambda F_0^{\gamma}}{\gamma} \Phi [F_0^{\gamma},1,1 + 
	  \lambda / \gamma] \right) = \mbox{constant}
\end{equation}
Consequently, there is no infinite tail in this situation, and the solution 
$X(F)$ intersects the $X$-axis, thus joining the trivial solution $F=0$. The 
reason for this lack of uniqueness at the boundary $F=0$ is that equation 
(\ref{prob3}) does not verify the Lipschitz condition on it, as it can be 
easily verified ({\sc Jackson}, 1991). Moreover, from (\ref{prob3}) we see that
\begin{equation}
\lim_{F \rightarrow 0^+} \dot{F}(F) = 
\lim_{F \rightarrow 0^+} \alpha F^g (1 - F^{\gamma} ) = 0
\end{equation}
and the solution is tangent to the $X$-axis at the intersection.

\subsubsection{Case III: $\;\: g=0 $}

This is a limit situation between the $g>0$ and the $g<0$ cases. Thus its 
interest is mainly mathematical, rather than applied. Since $\lambda =1$, 
we again have:
\begin{equation}
\label{noleft2}
	 \lim _{F \rightarrow 0^+} X(F) \; = \; X_0 -
	  \frac{F_0^{\lambda}}{\alpha \lambda} \left( 1 + 
	  \frac{\lambda F_0^{\gamma}}{\gamma} \Phi [F_0^{\gamma},1,1 + 
	  \lambda / \gamma] \right) = \mbox{constant}
\end{equation}
Therefore, there is not an infinite tail in this situation either, and the 
solution $X(F)$ again intersects the $X$-axis (but now $F=0$ is not a 
solution of the system, as it was in the case $0 < g < 1$). In addition, it 
is a simple task to check that the Lipschitz condition is not satisfied on 
the axis $F=0$ when $0<h<1$, thus indicating the existence of drawbacks for 
most numerical analyses of the solution in that region. Note also that now:
\begin{equation}
\lim_{F \rightarrow 0^+} \dot{F}(F) = 
\lim_{F \rightarrow 0^+} \alpha (1 - F^h ) = \alpha
\end{equation}
and the solution approaches the $X$-axis with slope $\alpha$.

\subsubsection{Case IV: $\;\: g<0 $}

Given that $\lambda >1$, we again have:
\begin{equation}
\label{noleft3}
	 \lim _{F \rightarrow 0^+} X(F) \; = \; X_0 -
	  \frac{F_0^{\lambda}}{\alpha \lambda} \left( 1 + 
	  \frac{\lambda F_0^{\gamma}}{\gamma} \Phi [F_0^{\gamma},1,1 + 
	  \lambda / \gamma] \right) = \mbox{constant}
\end{equation}
Therefore, an infinite left tail is present only in the generic Case I. As in 
Case III, now $F=0$ is not a solution of the system. However, the 
Lipschitz condition is never satisfied on the axis $F=0$ in this situation. 
Now we see that:
\begin{equation}
\lim_{F \rightarrow 0^+} \dot{F}(F) = 
\lim_{F \rightarrow 0^+} \alpha (F^g - F^h ) = + \infty
\end{equation}
independently of the value of $h$ (which now may be positive, zero or 
negative). Thus the solution joins the $X$-axis with infinite slope, i.e. 
perpendicularly.

In spite of their differences in what regards the left tail, it is interesting to mention 
that all the generic solutions (and the nongeneric ones as well, as we shall see) have an 
infinite right tail, independently of the values of the parameters. This is a direct 
consequence of the existence and uniqueness conditions verified by equation (\ref{prob2}) in 
the neighbourhood of the fixed point $F=1$.

Cases I, II and IV include most of the S-distributions that approximate
classical statistical distributions. Case I ($g>1$) includes Student
{\it t} distributions with few degrees of freedom and part of the region 
corresponding to non-central {\it t} distributions. Case II ($0<g<1$)
includes the approximation as S-distributions of many of the most common 
distributions: normal, F, $\chi^2$, non-central $\chi^2$, Weibull, $\chi^{-2}$, etc. 
Case IV is a secondary one, including some F-Distributions with (1,$\nu$) degrees of freedom,
and some restricted cases of non-central $\chi^2$ (see {\sc Voit, 1992}
for details). 

The four cases discussed so far include all practical situations
in data analysis since the nongeneric cases in which $g=1$ or 
$h=(1+k^{-1})g - k^{-1}, g>1$, and $k=(1,2,3...)$ do not arise in
data fitting, as we shall see in brief. Therefore, solution (\ref{gen2}) 
is the only one of interest in what concerns practical applications. 
However, {\it nongeneric} S-distributions, apart from having some interest
from the theoretical and the mathematical points of view, are indispensable 
in order to completely understand the generality of the {\it generic} case 
and the structure of solutions of the S-distribution. Such analysis is the 
aim of the next subsection.

\subsection{The nongeneric case}

\subsubsection{Structure of solutions in parameter space}

Before finding the explicit solution of the different nongeneric cases it is 
quite convenient to analyze how these are structured in parameter space. 
This is essential in order to clearly understand the structure of solutions 
of the S-distribution and leads to a natural structuration of them. Moreover, 
this is the clue for determining the relative importance of the generic and 
nongeneric cases in practical situations.

Let us recall that, for given $g$ and $h$ in (\ref{prob2}), we say that we
have a {\em nongeneric\/} or {\em degenerate\/} case if there is at least
one natural $k$ ($k=0,1,2 \ldots$) such that $k \gamma -g = -1$. In fact, 
such $k$ is always unique if exists. In other words, in the nongeneric case 
one and only one integral in the infinite series of integrals (\ref{intser}) 
produces a logarithmic term. 

We can rewrite the degeneracy condition as:
\begin{equation}
\label{deg}
	 \left\{ \begin{array}{lcl}
	 g = 1 \; , \:\; h > g & , &  k = 0 \\
	 h = \left( 1 + k^{-1} \right) g - k^{-1} \; , \:\; g > 1 
	 & , & k = 1, 2, 3 \ldots
	 \end{array} \right.
\end{equation}
These are, in all cases, straight lines that intersect at point $(1,1)$ 
of the $(g,h)$ plane, with infinite slope in the $k=0$ case, and slope 
$\:(1+1/k)\:$ in the $k$-th case (Figure 1). Notice also that $h>g$ and then 
the region of interest is one half of the plane $I \!\! R^2$.

\begin{center}
		  \mbox{} \\
		  \fbox{{\bf FIGURE 1 HERE (See Figure Captions List)}} \\
		  \mbox{} 
\end{center}

As a first remark, we see that nongeneric cases have zero measure when 
compared to the generic ones. Second, notice that the slopes of the 
half-lines corresponding to $k = 1, 2, 3 \ldots$ tend to 1 as $k 
\rightarrow \infty$. This means, as it can be seen in Figure 1, that the 
half-line 
\begin{equation}
\Omega \equiv \{ (g,h) \in I \!\! R^2 \;\: \mbox{such that} \;\: h = g  
\;\: \mbox{and} \;\:g > 1 \}
\end{equation}
is a limit set of the degenerate 
cases, which tend to $\Omega$ as $k \rightarrow \infty$ and are closely 
intermingled with the generic ones in every neighborhood of $\Omega$. 

This nice structure of the solutions poses an obvious difficulty when facing 
practical problems such as parameter estimation. However, such difficulty 
does not really arise in applications, since the distribution parameters resulting 
from data analysis will never lie on any of the zero-measure sets of degenerate 
solutions. Moreover, most distributions of  practical interest are usually far 
from every pathologic area ({\sc Voit}, 1992; {\sc Voit} and {\sc Yu}, 1994), since 
typical values are in the range $0<g<1.5$ and $g<h<8$, thus confirming our 
previous statement that solution (\ref{gen2}) is the only one of interest for 
practical purposes. 

This issue will be illustrated in the second part of this work. 
But before doing so we shall briefly consider the solutions of the 
degenerate problem. These can be divided in two qualitatively different 
possibilities:

\subsubsection{Case V (nongeneric): $\;\: g = 1 \;\: (k=0)$}

In this case we have $g=1$, $\lambda =0$, $h>1$ and $\gamma = h - 1 > 0$. From equation 
(\ref{intser}) we have:
\begin{equation} 
	 \sum_{k=0}^{\infty} \int_{F_0}^{F}  \xi^{k \gamma - g} \mbox{d} \xi =
	 \ln \left( \frac{F}{F_0} \right) + \sum_{k=1}^{\infty} 
	 \left. \frac{\xi^{k \gamma}}{k \gamma} \right|_{F_0}^{F}
\end{equation}
The right-hand side series can be summed, since 
\begin{equation} 
	\ln (1 + \eta) = \eta - \frac{1}{2} \eta^2 + \frac{1}{3} \eta^3 
	- \frac{1}{4} \eta^4 + \ldots  \;\:\:\;\:\; -1 < \eta \leq 1
\end{equation}
Thus the series is convergent in our case ($\eta = -F^{\gamma}$). Taking 
(\ref{sol1}) into account, we obtain:
\begin{equation} 
\label{solg1}
	X = X_0 + \frac{1}{\alpha} \left[ \ln \left( \frac{F}{F_0} \right) + 
	\frac{1}{\gamma} \ln \left( \frac{1-F_0^{\gamma}}{1-F^{\gamma}} \right)
	\right]
\end{equation}
This is the solution of the $g=1$ nongeneric case. Notice that now we will 
have an infinite tail in the left side of the distribution. In fact, the 
dominant behaviour is:
\begin{equation}
	 X (F \rightarrow 0^+) \; \sim \; \frac{\ln F}{\alpha} + \mbox{constant}
\end{equation}
Therefore $X \rightarrow - \infty$ when $F \rightarrow 0^+$, as announced.

\subsubsection{Case VI (nongeneric): $\;\: g \neq 1 \;\: (k=1,2,3 \ldots )$}

Now there exists a positive integer $N$, with $N \geq 1$, such that 
$N \gamma - g = -1$. We then find:
\begin{equation} 
	 \sum_{k=0}^{\infty} \int_{F_0}^{F}  \xi^{k \gamma - g} \mbox{d} \xi =
	 \sum_{k=0}^{N-1} \int_{F_0}^{F}  \xi^{k \gamma - g} \mbox{d} \xi 
	 + \int_{F_0}^{F}  \xi^{-1} \mbox{d} \xi
	 + \sum_{k=N+1}^{\infty} \int_{F_0}^{F}  \xi^{k \gamma - g} \mbox{d} \xi 
\end{equation}
Since this series is absolutely convergent in the domain of interest, we 
can regroup the terms and find, according to (\ref{sol1}):
\begin{equation} 
\label{solgk}
	X = X_0 + \frac{1}{\alpha} \left[ \ln \left( \frac{F}{F_0} \right) + 
	\sum_{k=0 \: , \: k \neq N}^{\infty} \frac{F^{k \gamma + \lambda} - 
	F_0^{k \gamma + \lambda}}{k \gamma + \lambda} \right]
\end{equation}
This is the solution of the $g \neq 1$ degenerate case. As in Case V, we have 
an infinite tail on the left of the distribution, as it can be inferred from 
the dominant behaviour:
\begin{equation}
	 X (F \rightarrow 0^+) \; \sim \; \frac{\ln F}{\alpha} + 
	 \frac{F^{\lambda}}{\alpha \lambda} + \mbox{constant}
\end{equation}
In other words, since $\lambda < 0$ in the present situation, we again have 
that $X \rightarrow - \infty$ when $F \rightarrow 0^+$.

Notice once again that in both nongeneric Cases V and VI, the distributions 
have infinite right tails for all parameter values.

This completes the classification and analysis of solutions of the 
S-distribution equation.

\section{Applications}

\subsection{Implementation of the quantile solution, computation of cumulatives and
design of special statistical distributions}

The S-distribution quantile solution as obtained in (\ref{gen2}), with the particular solutions
(\ref{solg1}) and (\ref{solgk}) for the nongeneric case, can be easily implemented. For instance,
we can use the program {\it Mathematica} since the Lerch's transcendent function is included in
the basic options (function {\bf LerchPhi}). As implemented in this software, this function 
may produce inadequate results for $v<0$ when computing (\ref{lerch}). According to (\ref{gen2}), 
in our case we have $v=1+\lambda/\gamma$, so the evaluation problems arise when $2 g>h+1$.
This is the region of the $(g,h)$ plane comprised between the $k=1$ and the $k=\infty$ 
degenerated solutions (see equation (\ref{deg})). In such cases, the correct result 
can be implemented by using the following property of Lerch's transcendent:

\begin{equation}
\label{PropiedadLerch}
\Phi[X^\gamma,1,v]=X^{\gamma m} \Phi[X^{\gamma},1,m+v] + \sum_{i=0}^{m-1} \frac{X^{\gamma i}}{v+i}
\end{equation}

\noindent For practical purposes a value of $m=100$ is appropriate. 
For a given value of the $cdf$ the corresponding quantile is immediately obtained
from (\ref{gen2}), (\ref{solg1}) or (\ref{solgk}), for a particular parameter set. 

Since the solution obtained is a quantile equation, the cumulatives shall be
computed by solving $X(F)=X_F$ for $F$. This can be easily done in {\it Mathematica} by
using {\bf FindRoot} and the quantile equation. However, for many practical purposes, the quantile
equation will suffice. For instance, given an S-distribution $S[F_0,X_0,\alpha,g,h]$ the 
{\it cdf} plot can be obtained by plotting the resulting quantiles from $F \in [0,1]$ and
reversing the resulting plot. 

Once the quantile solution is implemented, we can explore the theoretical results obtained in the
previous section. From a practical point of view, the most important result concerns to the
case $g<1$ in which no infinite left tail exists for the distribution. In such situations, 
the resulting S-distributions have a constrained range of possible values. The critical
value is $X(0)$, i.e. the value of $X$ at which the {\it cdf} attains 
a value of 0. For a specific parameter set in which $g<1$, it
is immediate to obtain the value of $X(0)$ evaluating the quantile equation with the value $F=0$.
For instance, in the case $S[0.5,10,1,.7,3]$ we obtain $X(0)=7.2211$. This critical
value for $X$ depends on the selected parameters and specially on $\alpha$, since this parameter is
inversely related with the variance.  

Although S-distributions with $g<1$ cannot, strictly speaking, reproduce exactly distributions with infinite 
left tails, they can model data with a long left tail provided that $\alpha$ has an appropriate value.
This is the case, for instance, when approaching the Normal distribution with an S-distribution 
since in this case $g \approx 0.69$ and $h \approx 2.88$ with $\alpha$ depending on $\sigma$ and $X(0.5)=\mu$
({\sc Voit}, 1992). 
The existence of a finite zero quantile in the case $g<1$ 
explains some problems that were observed when numerically computing the S-distribution solution using 
the differential equation (\ref{sddef}). For instance, using an S-distribution $S[0.5,10,1,.7,3]$ 
integration below $X=7.2211$ produces a numerical error. This problem has practical consequences for 
data analysis and makes the quantile analytical solution a useful alternative.

Far from being a problem, this property of the S-distributions with $g<1$ can be used to design truncated
distributions. Suppose we need a statistical distribution constrained to positive values of $X$, i.e. 
$X(0)=0$. To simplify the problem, we select $F_0=0.5$, so that $X_0$ will correspond
to the median of the distribution. Hence, we are left with four parameters: $X_0, \alpha, g$ and $h$.
To define the required distribution with $X(0)=0$, we first select three of the four S-distribution 
parameters and solve for the remaining one using the desired quantile value. For instance,
if we select $\alpha=1$, $g=0.1$ and $h=8$, the required initial condition with $F_0=0.5$ for having $X(0)=0$ 
is $X_0=0.595685$. The resulting S-distribution is shown in Fig. \ref{fig2}(a). A change in shape can be 
obtained by changing $g$ and $h$ and recomputing the initial condition. Following the same rationale, 
we can design a distribution that has a zero quantile for any specific value of $X$ different
from the origin. For instance, the example in Fig. \ref{fig2}(b) shows a distribution that 
is constrained to values greater than 20. Alternatively, we can choose a value for the median,
i.e.  the initial value $X_0$, fix two of the additional S-distribution parameters and
solve for the remaining one after imposing a given quantile value (Fig. \ref{fig2}(c)-(d)). This 
possibility can be very useful for using the resulting distributions in Monte-Carlo experiments. 
The extreme flexibility of the S-distributions makes these functions a valuable tool in 
this field. In the next section, we demonstrate the performance of the quantile solution 
in providing a way for generating random samples of a given S-distribution.

\begin{center}
		  \mbox{} \\
		  \fbox{{\bf FIGURE 2 HERE (See Figure Captions List)}} \\
		  \mbox{} 
\end{center}

\subsection{Random sample generation from an S-distribution}

Random sample generation of a given distribution requires using the inverse {\it cdf} equation, i.e.
the quantile equation. Since continuous evaluation of the differential equation may be very 
costly in terms of computer time, the search for alternatives is important to optimize quantile 
computations. {\sc Voit} and {\sc Schwacke} (2000) suggest using a simple rational function that very 
well approximates the quantile equation of the standard S-distribution. From this computation, 
that requires the use of tables containing coefficients for the different parameter values, 
other quantiles are easily computed. 

As an alternative, random data with the desired S-distribution behavior can be directly 
obtained by using the quantile solution represented by equations (\ref{gen2}), (\ref{solg1}) 
and (\ref{solgk}). Since the Lerch function is implemented in {\it Mathematica}, evaluation 
of the solution is quite efficient. In a standard Pentium II 300 MHz computer it takes about 
32 seconds to obtain a sample of 10000 S-distributed random numbers. In cases in which 
the solution is obtained by using (\ref{PropiedadLerch}), computer time may increase by
a factor of 10. However, this involves only those S-distributions with $2g>h+1$. The performance of this 
method can be seen in the examples included in Fig. \ref{fig3}. It may be appreciated that the 
samples obtained behave as expected by the theoretical distribution. In all cases, the samples 
are obtained in seconds and are ready for use in any procedure requiring data with a particular 
statistical distribution. 

\begin{center}
		  \mbox{} \\
		  \fbox{{\bf FIGURE 3 HERE (See Figure Captions List)}} \\
		  \mbox{} 
\end{center}

As stated in the previous subsection, we can artificially design an S-distribution so that one of the
quantiles is fixed. For instance, in Fig. \ref{fig4}(a) we defined a distribution that has only positive
values. This distribution can model a clinical parameter with values mostly within 0 and 1 and a
tail over 1. The generated random sample shows the expected behavior and can be used as simulated
data for exploring the effects of such statistical behavior on specific problems.  
In Fig. \ref{fig4} we present additional examples of special distributions designed to accomplish 
some specific features. For instance, in Fig. \ref{fig4}(b) we select a median of 50 and impose 
that the variable has only positive values. After selecting $g$ and $h$, we compute the required
$\alpha$ that results equal to 0.078126. Again, after selecting the required parameters, the
generated sample shows the expected statistical behavior.

\begin{center}
		  \mbox{} \\
		  \fbox{{\bf FIGURE 4 HERE (See Figure Captions List)}} \\
		  \mbox{} 
\end{center}

\subsection{Data representation using S-distributions}

Data representation using S-distributions has been discussed in detail elsewhere (see for 
instance {\sc Voit}, 1992, 2000; {\sc Balthis} et al. 1996; {\sc Sorribas} et al. 2000). In those papers, it has been shown that the S-distribution can be used as a parametric model for data despite the apparent complication of using a differential form that depends on four parameters. Parsimony arguments could be raised against this approach. However, the advantages related to a closed parametric form suggest that the S-distribution provides an useful solution. The results presented in the present paper help in solving some practical problems in the implementation of the S-distribution computations. Particularly, the quantile solution can be used as an altertative to the inverse S-distribution in those procedures that require computation of quantiles to fit an S-distribution to actual data. The two-step procedure suggested in {\sc Sorribas} et al. (2000) can be computationally improved by substituting the computations that require computation of the inverse S-distribution by the quantile solution. The final results are equivalent, but the required computer time is reduced by using the quantile solution.

We shall briefly recall the procedure and present some examples of the new implementation.
First, we organize the data in a histogram in which the total area equals one. The height of the
histogram at each class is taken as a $pdf$ value and the corresponding $cdf$ are computed
by adding the $pdf$ values. Taking $f_i=pdf_i$ and $F_i=cdf_i$ and using nonlinear regression 
we can fit:

\begin{equation}
\label{fitsd}
f_i=\alpha (F_i^g-F_i^h)
\end{equation}

\noindent and obtain a first estimation of $\alpha$, $g$ and $h$. Since these values are
taken as a first estimate, the choice of a bandwidth is not critical at this point.
Alternatively, the recently derived maximum likelihood estimator for the shape parameters $g$ and $h$ can be used 
to obtain an estimation of such parameters ({\sc Voit}, 2000). Once these parameters are computed, we refine the 
estimation by using a least-squares procedure that minimizes the 
sum of squares between the sample quantiles and the quantile values computed using the
corresponding S-distribution. Sample quantiles are the observed data points. The
corresponding sample {\it cdf} is used to compute the model quantiles using the quantile solution. 
The minimization procedure changes the S-distribution parameters until a satisfactory minimum is
reached. At this step, the previously estimated $g$ and $h$ are taken as fixed values and a new 
value for $\alpha$ and the initial condition $X_0$ are obtained. $F_0$ is taken as 0.5, so that
$X_0$ is an estimation of the median. This procedure provides
appropriate fits on simulated data (Fig. \ref{fig5}) and on actual data (Fig. \ref{fig6}). 
It is worth recalling here that S-distribution parameters have some built-in redundancy,
which leads to the existence of different parameter combinations that produce
S-distributions equivalent for practical purposes ({\sc Voit}, 1992; {\sc Sorribas} et al., 2000).
According to this property, the parameters obtained in Fig. \ref{fig6} once the fitting 
procedure is applied to a given data set may be different from the theoretical 
parameters of the original S-distribution. 

\begin{center}
		  \mbox{} \\
		  \fbox{{\bf FIGURE 5 HERE (See Figure Captions List)}} \\
		  \mbox{} 
\end{center}

\begin{center}
		  \mbox{} \\
		  \fbox{{\bf FIGURE 6 HERE (See Figure Captions List)}} \\
		  \mbox{} 
\end{center}

\section{Discussion}

The possibility of using a simple mathematical form as a general model for 
obtaining statistical distributions with any shape opens new possibilities for
data modeling and statistical simulations. The S-distribution provides such a
general model in terms of a differential equation on which the density is
a function of the cumulative. This particular definition leads to a simple
mathematical form but it had the drawback of requiring numerical integration 
for computing cumulatives and quantiles. 

In this contribution we have obtained the analytical solution for the quantile
equation of an S-distribution in terms of the Lerch function and it has been shown 
that this solution suffices for using this family of distributions in different 
practical applications. The analysis of the solution has led to identify some
interesting features of the S-distribution according to particular values of
parameter $g$. The most relevant observation concerns the case $g<1$ for which 
no infinite left tail exists. This result explains some problems that were observed 
when numerically integrating the S-distribution differential equation. Moreover, 
the possibility of computing the critical value $X_c$ that accomplishes 
$P(X \leq X_c)=0$ provides a way of designing S-distributions constrained to 
values over $X_c$. For instance, this may be the case of biological variables
that can have only posivite values. 

The quantile solution allows for identification of special cases in parameter space.
This is the case of a set of nongeneric solutions associated to $g=1$ and to values 
of $g>1$ for which $h=(1-k^{-1})g-k^{-1}$ with $k=1,2,3...$. In such situations,
we have derived the corresponding particular solutions and discussed their 
behavior. Altough additional results may be required for completely 
understanding the properties of the quantile solution obtained, our results 
clearly show the practical utility of this equation. On one hand, we can use it
to generate random samples for a given S-distribution. Such possiblity
allows the design of Monte-Carlo experiments in which the statistical
distributions can have almost any conceivable shape provided they are unimodal. Samples are obtained in
seconds with standard mathematical packages and can be used in applications. 
On the other hand, by using the quantile solution we can easily 
fit an S-distribution to observed data. The resulting fit can then be used
as a model for the population behavior of the observed variable.

The quantile solution presented in this paper should provide a tool for
optimizing some applications of the S-distribution in practical problems.
First, normality intervals for any variable are directly obtained
by computing the corresponding quantiles for the desired probability 
after fitting the corresponding S-distribution to the data set.
In addition, the simplicity in obtaining an S-distribution that fits
actual data using the quantile solution can facilitate the charaterization
of trends on the distribution of a given property as a function of 
some variable of interest. This may be the case of deriving weight 
and height standard curves as a function of age, for instance.

The concept of S-distribution opens a new way of representing a random 
variable in which a single mathematical form suffices for any particular case. 
To attain this flexibility, the price of using a differential equation that 
complicates its practical use is to be paid. Although this may be seen as a 
problem from a theoretical point of view, it is an advantage when it comes 
to distribution design and data modeling. The need of numerical methods for
using the S-distribution is now highly simplified by the analytical solution
obtained in this paper. We hope this can facilitate the obtainment of new
results and the application of the S-distribution to practical problems.

\vspace{1cm}

\noindent{\bf Acknowledgements}
\noindent{\em The authors are grateful to Dr.E.O.Voit for his criticism of the
original manuscript and for numerous suggestions that improved its final version. 
A.S. acknowledges support from the FIS of Spain (grant 00/0235) and a grant from 
La Paeria (Lleida X0148). B. H. acknowledges financial support from the 
University of Lleida.}

\mbox{}

\mbox{}

\pagebreak
\noindent{\bf References}

\mbox{}

\noindent {\sc Balthis, W. L.}, {\sc Voit, E. O.} and {\sc Meaburn, G. M.}, 1996:
Setting prediction limits for mercury concentrations in fish having high
bioaccumulation potential. {\em Environmetrics\/} {\bf 7}:429-439.
\vspace{3mm}

\noindent {\sc Churchill, R. V.} and {\sc Brown, J. W.}, 1984: 
{\em Complex Variables and Applications,\/} Fourth Edition. McGraw-Hill, 
New York, Chapter 5.
\vspace{3mm}

\noindent {\sc Erd\'{e}lyi, A., Magnus, W., Oberhettinger, F.} and {\sc Tricomi, F. G.}, 1953: 
{\em Higher Transcendental Functions,\/} Vol. I. McGraw-Hill, New York.
\vspace{3mm}

\noindent {\sc Jackson, E. A.}, 1991: 
{\em Perspectives of nonlinear dynamics\/}, Vol. I. Cambridge University Press, 
Cambridge, UK, Section 2.2.
\vspace{3mm}

\noindent {\sc Johnson, N. L.} and {\sc Kotz, S.}, 1970: 
{\em Continuous univariate distributions-1.\/} Houghton Mifflin, Boston, MA.
\vspace{3mm}

\noindent {\sc Ryu, H.K.}, 1993: 
Maximum entropy estimation of density and regression functions. 
{\em J.Econometrics.\/} {\bf 3}:379-440.
\vspace{3mm}

\noindent {\sc Sorribas, A., March, J.} and {\sc Voit, E. O.}, 2000: 
Estimating age-related trends in cross-sectional studies using S-distributions. 
{\em Stat. Med.\/} {\bf 19}:697-713.
\vspace{3mm}

\noindent {\sc Voit, E. O.}, 1992: 
The S-distribution: A tool for approximation and classification of univariate, 
unimodal probability distributions. {\em Biom. J.\/} {\bf 7}:855-878.
\vspace{3mm}

\noindent {\sc Voit, E. O.}, 2000: A Maximum Likelihood Estimator for Shape Parameters of 
S-distributions. {\em Biom. J.\/} {\bf 42}:471-479.
\vspace{3mm}

\noindent {\sc Voit, E. O., Balthis, W. L.} and {\sc Holser, R. A.}, 1995: Hierarchical 
Monte Carlo modeling with S-distributions: concepts and illustrative analysis of 
mercury contamination in king mackerel. {\em Environment Int.\/} {\bf21}:627-635.
\vspace{3mm}

\noindent {\sc Voit, E. O.} and {\sc Savageau, M. A.}, 1984: 
Analytical solution to a generalized growth equation. {\em J. Mathem. Anal. Appl.\/} 
{\bf103}:380-386.
\vspace{3mm}

\noindent {\sc Voit, E. O.} and {\sc Schwacke, L.H.}, 2000: 
Random Number Generation from Right-Skewed, Symmetric,
and Left-Skewed Distributions. {\em Risk Analysis.\/} 
{\bf 20}:393-402.
\vspace{3mm}

\noindent {\sc Voit, E. O.} and {\sc Sorribas, A.}, 2000: 
Computer modeling of dynamically changing distributions of random variables. 
{\em Math. Comp. Modeling\/} {\bf 31}:217-225.
\vspace{3mm}

\noindent {\sc Voit, E. O.} and {\sc Yu, S.}, 1994: 
The S-Distribution: Approximation of Discrete Distributions. 
{\em Biom. J.\/} {\bf 36}: 205--219.

\noindent {\sc Wagner, U.} and {\sc Geyer, A.L.J.}, 1995: 
A maximum entropy method for inverting Laplace transforms of probability density functions. 
{\em Biometrika.\/} {\bf 4}:887-892.
\vspace{3mm}

\noindent {\sc Wang, M.P.} and {\sc Jones, M.C.}, 1995: 
{\em Kernel smoothing,\/}  Chapman and Hall, London.
\vspace{3mm}

\pagebreak
\thispagestyle{empty}
\begin{center}
	 {\large {\bf FIGURE CAPTIONS}}
\end{center}

\mbox{}

\begin{fig} 
\label{fig1} {\rm Structure of solutions of the S-distribution equation on 
the $g-h$ plane. The continuous thick black lines correspond to degeneracies 
(numbered in terms of $k$, see equation (\ref{deg})). The different regions (I-VI) 
classified in the main text are also indicated. Region III is the $h$-axis 
(black dashed line). The lower half plane $h \leq g$ is excluded by 
definition (in grey). 
}
\end{fig}

\begin{fig} 
\label{fig2} {\rm Design of special statistical distributions. Using the 
quantile solution, we can design statistical distributions that fulfill 
some constraints. In all cases we take $F_0 = 0.5$ and represent the 
corresponding $pdf$.
\newline
(a) We select $\alpha=1$, $g=0.1$ and $h=8$ and
impose $X(F=0)=0$, i.e. the variable has only positive values. Solving the 
quantile equation, the required initial value is $X_0=0.595685$. 
\newline
(b) We select $g=0.1$, $h=8$ and an initial condition $X_0=50$. If we impose
$X(F=0)=20$ and solve the quantile equation for $\alpha$ we obtain a value
of $\alpha=0.0198562$. 
\newline
(c) We select $\alpha=0.5$, $h=3$ and an initial 
condition $X_0=50$. If we impose $X(F=0.1)=12$ and solve the quantile 
equation for $g$ we obtain a value of $g=2.28146$. 
\newline
(d) We select 
$\alpha=0.5$, $h=3$ and an initial condition $X_0=50$. If we impose 
$X(F=0.1)=45$ and solve the quantile equation for $g$ we obtain a value of 
$g=1.2235$. 
}
\end{fig}

\begin{fig} 
\label{fig3} {\rm Random sample generation from an S-distribution.
Samples are obtained by generating a random sample from an uniform
distribution and by computing the corresponding quantiles for the 
selected S-distribution. (a) $S[0.5, 100, 0.2, 1, 3]$; (b) 
$S[0.5, 100, 0.2, 1, 7]$; (c) $S[0.5, 100, 0.2, 0.1, 7]$; 
(d) $S[0.5, 100, 0.1, 0.4, 7]$. In each case, 
we plot the histogram of the data set and the corresponding $pdf$ for
the selected S-distribution. Histograms are scaled so that the total
area adds to one.
}
\end{fig}

\begin{fig} 
\label{fig4} {\rm  Random sample generation from S-distribution designed
to fulfill special requirements. Data are obtained following the same
procedure indicated in Fig. \ref{fig3}. \newline
(a) $S[0.5, 0.595685, 1, 0.1, 0.8]$. Parameters are choosen so that $X(F=0)=0$. 
\newline
(b) $S[0.5, 36.8252, 0.1, 0.69, 2.88]$. Parameters are choosen so that $X(F=0)=10$. 
\newline
(c) $S[0.5, 8.90629, 0.1, 0.3, 4]$. Parameters are choosen so that $X(F=0)=0$. 
\newline
(d) $S[0.5, 50, 0.0178126, 0.3, 4]$. Parameters are choosen so that $X(F=0)=0$.
\newline
}
\end{fig}

\pagebreak
\begin{fig} 
\label{fig5} {\rm Fitting an S-distribution to observed data. Data are generated
from an S-dis\-tri\-bu\-tion following the same procedure indicated in Fig. \ref{fig3}.
Then, an S-distribution is fitted using the procedure indicated in the text.
In each case, the grey line indicates the theoretical distribution and the black
line indicates the fitted distribution. \newline
(a) Theoretical distribution: $S[0.5, 50, 1, 0.6, 3]$. Sample size $= 200$.
Fitted distribution: $S[0.5, 50.0409, 1.0491, 0.694497, 3.40337]$.
\newline
(b) Theoretical distribution: $S[0.5, 50, 1, 0.2, 3]$. Sample size $= 200$.
Fitted distribution: $S[0.5, 50.0298, 1.09488, 0.356604, 3.32489]$.
\newline
(c) Theoretical distribution: $S[0.5, 50, 1, 0.6, 7]$. Sample size $= 500$.
Fitted distribution: $S[0.5, 50.0203, 0.995178, 0.61789, 9.44785]$.
\newline
(d) Theoretical distribution: $S[0.5, 100, 0.2, 0.1, 7]$. Sample size $= 300$.
Fitted distribution: $S[0.5, 100.006, 0.216268, 0.161938, 6.24436].$
}
\end{fig}

\begin{fig} 
\label{fig6} {\rm Fitting an S-distribution to actual clinical data. 
Data were collected by the Intensive Care Service of the Hospital
Arnau de Vilanova in Lleida (Spain) between 1996 and 1998. 
All patients (n=507) had been in intensive care more than 24 hours. 
The parameters of the fitted S-distributions are:
\newline
(a) $S[0.5, 53.4528, 1.117, 0.7158, 0.7386]$;
\newline
(b) $S[0.5, 3.83737, 21.106, 1.04544, 1.11902]$; 
\newline
(c) $S[0.5, 7.36712, 8.25116, 0.944545, 4.44792]$;
\newline
(d) $S[0.5, 98.7516, 0.632082, 1.1001, 1.14905]$. 
}
\end{fig}

\pagebreak

Author's address:

Albert Sorribas

Departament de Ci\`{e}ncies M\`{e}diques B\`{a}siques

Universitat de Lleida

Av.Rovira Roure, 44

25198-Lleida

Spain

e-mail: Albert.Sorribas@cmb.udl.es
\end{document}